\documentclass{gtpart}
\usepackage{amsmath,amsthm,amsfonts,amssymb,amscd,enumerate}

\title[Examples of buildings]{Examples of buildings constructed via covering spaces}
\author[M\,W Davis]{Michael W. Davis}
\givenname{Michael}
\surname{Davis}
\address{The Ohio State University\\
Department of Mathematics\\
231 W. 18th Ave.\\
Columbus, OH 43210-1174\\
United States}
\email{mdavis@math.ohio-state.edu}
\urladdr{http://www.math.ohio-state.edu/~mdavis/}

\keyword{Building,  Coxeter group}
\subject{primary}{msc2000}{20F65}
\subject{primary}{msc2000}{20F55}
\subject{primary}{msc2000}{20E42}
\subject{secondary}{msc2000}{20J06}
\subject{secondary}{msc2000}{51E24}
\subject{secondary}{msc2000}{57M07}
\arxivreference{0809.5097}
\arxivpassword{w7c7u}

\volumenumber{}
\issuenumber{}
\publicationyear{}
\papernumber{}
\startpage{}
\endpage{}
\doi{}
\MR{}
\Zbl{}
\received{}
\revised{}
\accepted{}
\published{}
\publishedonline{}
\proposed{}
\seconded{}
\corresponding{}
\editor{}
\version{}

\theoremstyle{plain}
\newtheorem{theorem}{Theorem}[section]
\newtheorem{proposition}[theorem]{Proposition}

\newtheorem{corollary}[theorem]{Corollary}
\newtheorem{lemma}[theorem]{Lemma}

\newtheorem*{Question}{Question}

\theoremstyle{definition}

\newtheorem{definition}[theorem]{Definition}

\newtheorem{example}[theorem]{Example}
\newtheorem{examples}[theorem]{Examples}

\newtheorem*{Remark}{Remark}
\newtheorem*{remarks*}{Remarks}

\newcommand{\ca}{\mathcal {A}}

\newcommand{\cac}{\mathcal {C}}

\newcommand{\car}{\mathcal {R}}
\newcommand{\cs}{\mathcal {S}}
\newcommand{\ct}{\mathcal {T}}
\newcommand{\cu}{\mathcal {U}}

\newcommand{\ff}{{\mathbb F}}

\newcommand{\rr}{{\mathbf R}}

\newcommand{\bprime}{{\mathbf s'}}
\newcommand{\zz}{{\mathbf Z}}

\newcommand{\minus}{{-1}}

\newcommand{\K}{{C_{\langle K \rangle}}}

\newcommand{\bC}{{\mathbf C}}

\newcommand{\bs}{{\mathbf s}}

\newcommand{\tcu}{{\widetilde{\cu}}}

\newcommand{\wt}{\widetilde}

\newcommand{\ol}{\overline}

\newcommand{\gd}{\delta}

\newcommand{\gf}{\varphi}

\newcommand{\gr}{\rho}
\newcommand{\gs}{\sigma}

\newcommand{\gG}{\Gamma}
\newcommand{\gD}{\Delta}

\newcommand{\gO}{\Omega}

\newcommand{\Aut}{\operatorname{Aut}}
\newcommand{\Card}{\operatorname{Card}}
\newcommand{\cent}{\operatorname{Cent}}

\newcommand{\edge}{\operatorname{Edge}}

\newcommand{\Ker}{\operatorname{Ker}}

\newcommand{\flag}{\operatorname{Flag}}

\newcommand{\vertex}{\operatorname{Vertex}}
\newcommand{\Prod}{\operatorname{prod}}

        {\end{enumerate}}
        {\end{enumerate}}       
\newenvironment{enumerate1}{\begin{enumerate}[\upshape 1)]}%
        {\end{enumerate}}
\newenvironment{enumerateI}{\begin{enumerate}[\upshape 
(I)]}{\end{enumerate}}

\newenvironment{enumeratea'}{\begin{enumerate}[\upshape 
(a)$'$]}{\end{enumerate}}

\newcommand{\comment}[1]{}

\def\l{\operatorname{\ell}}

\providecommand{\bysame}{\makebox[3em]{\hrulefill}\thinspace,\ }

\numberwithin{equation}{section}
\begin{document}



\begin{abstract}
Covering space theory is used to construct new examples of buildings.
\end{abstract}

\maketitle

\section*{Introduction}
A building consists of a Coxeter system $(W,S)$, a set $\cac$ (of ``chambers'') and a ``Weyl distance'' $\gd:\cac\times \cac \to W$, satisfying certain axioms.  Given a building $\cac$,a space $X$ and a ``mirror structre'' on $X$ , one can construct a ``realization'' $\cu(\cac,X)$ of the building, by pasting together copies of $X$, one for each element of $\cac$.  (Details are given in \S \ref{ss:b}.)  The classical choice for $X$ is a simplex.  

The following types of buildings have been studied:
\begin{itemize}
\item
spherical buildings, where $W$ is a finite Coxeter group,
\item
affine buildings, where $W$ is a Euclidean reflection group,
\item
Kac-Moody buildings over finite fields, where $W$ is fairly arbitrary but where the ``thickness''  is restricted to be a constant prime power $q$,
\item
right-angled buildings, where $W$ is required to be a right-angled Coxeter group.
\end{itemize}

The purpose of this note is to use covering space theory to construct  new examples of buildings.  On the level of Coxeter systems, this contsruction was described previously in Davis \cite{d87}.  It goes as follows.  Start with a Coxeter system $(W',S)$.  Change some $m(s,t)$'s in its Coxeter matrix from integers $\ge 2$ to the symbol $\infty$.  This defines a new Coxeter system $(W,S)$.  Suppose $\cac'$ is a building of type $(W',S)$.  We want to use it construct a building $\cac$ of type $(W,S)$.  Let $K$ be the ``Davis chamber'' for $(W, S)$.  Use it to get a realization $\cu(\cac',K)$ for $\cac'$.  Since we made changes in $W'$, the space $\cu(\cac', K)$ will not be simply connected.  However, as we show in Theorem~\ref{t:cac}, its universal cover is the realization of a building $\cac$ of type $(W,S)$. 

Interesting examples can be constructed by starting with the case where $W'$ is a finite Coxeter group and $\cac'$ is a spherical building. In particular, there are examples of buildings with the following features. 
\begin{itemize}
\item
There are buildings whose maximal spherical residues contain an arbitrary set of spherical buildings and whose standard realiztions are highly connected at infinity (ruling out the procedure of simply taking the free product of various spherical buildings).
\item
In contrast to the (irreducible) Kac-Moody case, there are examples of buildings with different thicknesses for different types  of  spherical residues (in other words, even when these spherical residues are associated to algebraic groups over finite fields, the fields may be different).
\item
When $\cac'$ is a finite spherical building, the automorphism group, $\Aut (\cac)$ always contains a torsion-free uniform lattice. (Again, this is in contrast to the Kac-Moody case: Kac-Moody groups generally are not expected to contain any uniform lattice.)
\end{itemize} 
In \S 3 we discuss the special case where $\cac'$ is a direct product, $\cac_0\times\cdots \times \cac_p$. 
By making appropriate choices of which $m(s,t)$'s to change to $\infty$, we can define the notions of  ``free products'' and ``graph products'' of buildings analogously to the corresponding notions for groups.  When each of the buildings $\cac_i$ has the same rank, there also is the  notion of a ``square product,'' $\cac_0\, \square \cdots\square \,\cac_p$ (see Example~\ref{ex:square}).  The square product is interesting even when there are only two factors, $\cac_0$, $\cac_1$, and both are spherical buildings.  In this case the fundamental chamber $K$ is combinatorially isomorphic to an $n$-cube.  The link of one vertex of the cube in the standard realization, $\cu(\cac_0\, \square\,\cac_1,K)$, is the spherical realization of $\cac_0$ while the link at an opposite vertex is the spherical realization of $\cac_1$ (see Example~\ref{ex:cube}). 

\emph{Acknowledgements}\ \ My thanks go to Tadeusz Januszkiewicz for several helpful conversations concerning these examples over the past twenty years.  While writing this paper I was  partially supported by NSF 
grant DMS 0706259.

\section{Basic definitions}

\subsection{Coxeter groups}\label{ss:cg}
A \emph{Coxeter matrix over a set} $S$ is an $S\times S$ symmetric matrix $M=(m(s,t))$ with each diagonal entry $1
$ and each off-diagonal entry either an integer $\ge 2$ or the symbol $\infty$.  The matrix $M$ defines a presentation of a group $W$ as follows:  the set of generators is $S$ and the relations have the form $(st)^{m(s,t)}$ where $(s,t)$ ranges over all pairs in $S\times S$ such that $m(s,t)\neq \infty$.  The pair $(W,S)$ is a \emph{Coxeter system} (cf Bourbaki \cite{bourbaki}).
Given  $T\subset S$, $W_T$ denotes the subgroup generated by $T$; it is called a \emph{special subgroup}.  $(W_T,T)$ is itself a Coxeter system (see Bourbaki \cite[IV.8, Thm.~2(i)]{bourbaki}).  The subset $T$ is \emph{spherical} if $W_T$ is finite.  

\begin{definition}
The poset of spherical subsets of $S$ (partially ordered by inclusion) is denoted $\cs$ (or $\cs(W,S)$).  The set of nonempty elements in $\cs$ is an abstract simplicial complex $L$ (or $L(W,S)$) called the \emph{nerve} of $(W,S)$.  (In other words, the vertex set of $L(W,S)$ is $S$ and a nonempty subset $T\subset S$ spans a simplex if and only if it is spherical.)
\end{definition}

Given a word $\bs=(s_1,\dots, s_n)$ in $S$, its \emph{value}, $w(\bs)$, is the element of $W$ defined by $w(\bs):=s_1\cdots s_n$.  The word $\bs$ is a \emph{reduced expression} if it is a word of minimum length for $w(\bs)$, ie, if $l(w(\bs))=n$.

\textbf{Tits' solution to the word problem}\ \  
Given elements $s,t\in S$ and an integer $m\geq 2$, let  $\Prod (s,t;m)$ denote the alternating word $(s,t,\dots)$ in $s$ and $t$ beginning with $s$ and having length $m$.  Consider the following two operations on a word $\bs$:
\begin{enumerateI}
\item
Delete a subword of the form $ss$
\item
Replace a subword $\Prod (s,t;m(s,t))$ by $\Prod (t,s;m(s,t))$.
\end{enumerateI} 
An operation of type (I) is a \emph{cancellation}; one of type (II) is a \emph{flip}.

Tits proved that if a word $\bs=(s_1,\dots, s_n)$ is not a reduced expression for an element $w$, then it can be changed to a reduced expression for the same element by a sequence of flips and cancellations.  Moreover, if $\bs$ and $\bprime$ are two reduced expressions for the same element, then one can be changed to the other by a sequence of flips (see Abramenko--Brown \cite[\S2.3.3]{ab}).

\subsection{Buildings}\label{ss:b}

\begin{definition}\label{d:bldg}
(cf Abramenko--Brown \cite{ab} or Ronan \cite{ronan}). Suppose $(W,S)$ is a Coxeter system. 
A \emph{building of type $(W,S)$} is a pair $(\cac, \gd)$ consisting of a nonempty set $\cac$ (the elements of which are called \emph{chambers}), and a function $\gd:\cac\times\cac\to W$ (called the \emph{Weyl distance}) so that the following conditions hold for all $C,D\in \cac$

\begin{enumerate}
\item[]
\textbf{(WD\upshape 1)} 
$\gd(C,D)=1$ if and only if $C=D$.
\item[]
\textbf{(WD\upshape 2)} 
If $\gd(C,D)=w$ and $C'\in\cac$ satisfies $\gd(C',C)=s\in S$, then $\gd(C',D)=sw$ or $w$.  If, in addition, $\l(sw)=l(w)+1$, then $\gd(C', D)=sw$.
\item[]
\textbf{(WD\upshape 3)} If $\gd(C,D)=w$, then for any $s\in S$ there is a chamber $C'\in \cac$ such that $\gd(C',C)=s$ and $\gd(C', D)=sw$.
\end{enumerate}
\end{definition}

\begin{example}\label{ex:thin}
The group $W$ itself has the structure of a building:  $\gd:W\times W\to W$ is defined by $\gd(v,w)=v^\minus w$.  This is called the \emph{thin building} of type $(W,S)$.
\end{example}
A \emph{chamber system} is a set $\cac$ together with a family of equivalence relations on $\cac$ indexed by another set $S$.  Chambers $C,D\in \cac$ are \emph{$s$-adjacent} if they are $s$-equivalent and not equal.
A \emph{gallery} in $\cac$ is a finite sequence of 
chambers $(C_0, \dots ,C_k)$ such that $C_{j-1}$ is adjacent to $C_j, 1 \le j \le k$.  The  \emph{type} of this gallery is the word $\bs
=(s_1,  \dots, s_k)$ where $C_{j-1}$ is $s_j$-adjacent to $C_j$.  If each $s_j$ 
belongs to a given subset $T$ of $S$, then the gallery is a \emph{$T$-gallery}. Two chambers are in the \emph{same $T$-connected component} if they can be connected by a $T$-gallery.
The $T$-connected components of a chamber system $\cac$ are its 
\emph{residues} of  \emph{type $T$}.  An $s$-equivalence class is the same thing as a residue of type $\{s\}$. (An $s$-equivalence class is sometimes called a \emph{panel}.)

Now suppose $\cac$ is a building of type $(W,S)$.
Call chambers $C,D\in \cac$ \emph{$s$-adjacent} if $\gd(C,D)=s$; they are \emph{$s$-equivalent} if they are either $s$-adjacent or equal. This gives $\cac$ the structure of a chamber system.
When $\cac=W$, a $T$-residue is just a left coset of $W_T$.

A residue in a building of type $T$ is itself a building of type  $(W_T,T)$.  
A building of type $(W,S)$ is \emph{spherical} if $W$ is finite.  A building has \emph{finite thickness} if each $s$-equivalence class is finite, for each $s\in S$.  (This implies all spherical residues are finite.)  \emph{Henceforth, all buildings will be assumed to have finite thickness}.

A gallery $(C_0,\dots , C_k)$ of type $\bs$ is \emph{reduced} if $\bs$ is a reduced expression for $w(\bs)$.

It is proved in Abramenko--Brown  \cite[Prop.~5.23]{ab} that conditions \textbf{(WD1)},\textbf{(WD2)}, \textbf{(WD3)} in Definition~\ref{d:bldg} are equivalent to the following two conditions on a chamber system $\cac$ equipped with a function $\cac\times\cac\to W$.
\begin{itemize}
\item
Each $s$-equivalence class has at least two elements.
\item
Given a reduced expression $\bs$ for an element $w\in W$,  there is a gallery of type $\bs$ from $C$ to $D$ if and only if $\gd(C,D)=w$. 
\end{itemize}
(This is the definition in Ronan \cite{ronan}.)

An \emph{apartment} in $\cac$ is a subset which is $W$-isometric to the thin building $W$.  In other words, it is a subset $\ca\subset \cac$ such that for any $C\in \ca$, the function $\gr_C:\ca\to W$ defined by $D\mapsto \gd(C,D)$ is an isomorphism.

An \emph{automorphism} of $\cac$ is a self-bijection which preserves Weyl distance.  Given $C\in \cac$ the \emph{combinatorial ball of radius $n$ about $C$} is the set 
\[
B_C(n):=\{D\in \cac\mid l(\gd(C,D))\le n\}.  
\]
There is a natural topology on the group $\Aut (\cac)$ of automorphisms of $\cac$:  a small open neighborhood of $1\in \Aut (\cac)$ is the set of automorphisms which fix each element of $B_C(n)$ (for large $n$ and some $C$).  Since $\cac$ is locally finite, $\Aut(\cac)$ is a locally compact, totally disconnected topological group.  As such, it has a Haar measure.  A closed subgroup $G\subset \Aut(\cac)$ inherits a topology and a Haar measure.  A subgroup $\gG\subset G$ is a \emph{lattice} if it is discrete and $G/\gG$ has finite volume.  It is a \emph{uniform} lattice in $G$ if $G/\gG$ is compact.  A discrete subgroup $\gG\subset \Aut (\cac)$ is a uniform lattice if and only if $\cac/\gG$ is finite.  (See Thomas \cite{thomas} for a discussion of lattices in $\Aut (\cac)$, when $\cac$ is a right angled building.)

\begin{definition}
A subgroup $G\subset \Aut (\cac)$ is \emph{chamber-transitive} if is transitive  on $\cac$.  It is \emph{strongly transitive} if it is transitive on the set of pairs $(\ca, C)$, where $\ca$ is an apartment in $\cac$ and $C\in \ca$ (cf Abramenko--Brown \cite[\S 6.1.1]{ab}).  (In fact, it is not necessary use all apartments in this definition, $\ca$ need only belong to a certain ``system of apartments'' saisfying the classical axioms for a building, cf \cite[\S6.1]{ab}).
\end{definition}

It turns out that if $G$ is strongly transitive on a thick building, then it inherits the structure of a $BN$ pair (also called a ``Tits system''), cf \cite[Thm.~6.56]{ab}. (A building is \emph{thick} if each $s$-equivalence class has at least 3 elements.)

\subsection{Geometric realizations}\label{ss:gr}

A \emph{mirror structure} over a set $S$ on a space $X$ is a family of subspaces $(X_s)_{s\in S}$ indexed by $S$.  Given a mirror structure on $X$, a subspace $Y\subset X$ inherits a mirror structure by setting $Y_s:=Y\cap X_s$.  If $X$ is a CW complex and each $X_s$ is a subcomplex, then $X$ is a \emph{mirrored CW complex}.  For each nonempty subset $T\subset S$, define subspaces $X_T$ and $X^T$ by
\begin{equation}\label{e:T}
X_T:=\bigcap_{s\in T} X_s\quad\text{and}\quad X^T:=\bigcup_{s\in T} X_s. 
\end{equation}
Put $X_\emptyset:=X$ and $X^\emptyset:=\emptyset$.  
Given a cell $c$ of (a CW complex) $X$ or a point $x \in X$, put
\begin{align*}
S(c)&:=\{s\in S\mid c\subset X_s\},\\
S(x)&:=\{s\in S\mid x\in X_s\}.
\end{align*}

Given a building $\cac$ of type $(W,S)$ and a mirrored space $X$ over $S$,
define an equivalence relation  $\sim$ on $\cac \times X$ by $(C, x) \sim (D, y)$ if and 
only if $x = y$ and $\gd (C, D) \in W_{S(x)}$ (ie, $C$ and $D$ belong to the same $S(x)$-residue).  The $X$-\emph{realization} 
of $\cac$, denoted  $\cu(\cac,X)$, is defined by
	\begin{equation}\label{e: defUPhi}
	\cu(\cac,X):=(\cac \times X)/\sim.  
	\end{equation}
($\cac$ has the discrete topology.)  

\begin{example}\label{ex:spherical}
(\emph{Spherical realizations}) 
Suppose $W$ is a finite Coxeter group of rank $n+1$ and $\cac$ is a spherical building of type $(W,S)$.  Then $W$ is a reflection group on the unit sphere $S^n\subset \rr^S$.  Let $\gD^n\subset S^n$ be the fundamental simplex with its codimension one faces indexed by $S$. $\cu(\cac,\gD^n)$ is called the \emph{spherical realization} of $\cac$.  It is homotopy equivalent to a wedge of $n$-spheres (cf Abramenko--Brown \cite[Thm.~4.127]{ab}).
\end{example}

\textbf{Geometric realizations of posets}\ \   
Given a poset $\ct$, $\flag(\ct)$ denotes the set of finite chains in $\ct$, partially ordered by inclusion, ie, an element of $\flag(\ct)$ is a finite, nonempty, totally ordered subset of $\ct$.  
$\flag(\ct)$ is an abstract simplicial complex with vertex set $\ct$ and with $k$-simplices the elements of $\flag(\ct)$ of cardinality $k+1$.  The corresponding topological simplicial complex is the \emph{geometric realization} of the poset $\ct$ and is denoted by $|\ct|$.  

Suppose $L$ is a simplicial complex with vertex set $S$.  Let $\cs(L)$ be the poset of subsets of $S$ (including $\emptyset$) which are vertex sets of simplices in $L$.  Let $K(L):=|\cs(L)|$ be the geometric realization of this poset.  ($K(L)$ is isomorphic to the cone on the barycentric subdivision of $L$, where $\emptyset\in \cs$ provides the cone point.) Define a mirror structure on $K(L)$ by putting $K(L)_s:=|\cs(L)_{\geq \{s\}}|$.

\textbf{Galleries in $\boldsymbol{\cac}$ and edge paths in $\boldsymbol{\cu(\cac,K(L))}$}\ \  Let $v_\emptyset$ denote  the central vertex of $K(L)$ corresponding to $\emptyset$ and for each $C\in\cac$, let $v_C$ denote the image of $(C,v_\emptyset)$ in $\cu(\cac,K(L))$.  The vertices $v_C$, with $C\in \cac$, are the \emph{central vertices} of $\cu(\cac,K(L))$.  Let $v_s$ denote the vertex of $K(L)$ coresponding to $\{s\}$.  If $\car$ is a residue of type $\{s\}$ (ie a ``panel''), then for any $C\in \car$ the image of $(C,v_s)$ in $\cu(\cac,K(L))$ is denoted $v_\car$ and called a \emph{mirror vertex}.  The union of edges in $\cu(\cac,K(L))$ connecting central vertices to mirror vertices is called the \emph{dual $1$-skeleton} of $\cu(\cac,K(L))$.  Given a gallery $(C_0,\dots, C_k)$ in $\cac$ we get an edge path in the dual $1$-skeleton with successive vertices $(v_{C_0},v_{\car_1},v_{C_1},\dots, v_{C_k})$.  (There are $2k+1$ vertices in this edge path; every other one has the form $v_{\car_j}$ where $\car_j$ is a residue of type $s_j$,)  We will simplify notation by omitting the odd vertices and writing $(v_{C_0},\dots, v_{C_k})$ for the corresponding edge path in $\cu(\cac,K(L))$.

\textbf{The standard realization, $\boldsymbol{\cu(\cac,K)}$}\ \ 
As before, $\cs$ denotes the poset of spherical subsets of $S$. Put $K:=|\cs|$.   For each $s\in S$, put $K_s:=|\cs_{\ge \{s\}}|$ and for each $T\in \cs$, $K_T=|\cs_{\ge T}|$.  $K$ is sometimes called the \emph{Davis chamber} of $(W,S)$ and $\cu(W,K)$, the \emph{Davis complex}.  The simplicial complex $\cu(\cac, K)$ is the \emph{standard realization} of $\cac$.  (Alternatively, $\cu(\cac,K)$ is the geometric realization of the poset of spherical residues of $\cac$ (see Davis \cite{dbuild}).
By construction $\cu(\cac, K)$ is locally finite (since $\cac$ is assumed to have finite thickness).  It is proved in \cite{dbuild} that $\cu(\cac,K)$ is contractible. 

\subsection{Duality groups}\label{ss:duality}
A space $X$ is an \emph{$n$-dimensional duality space} if it is acyclic and if its compactly supported cohomology, $H^*_c(X)$, is torsion-free and concentrated in degree $n$. A discrete group $\gG$ is an \emph{$n$-dimensional duality group} if $H^*(\gG;\zz\gG)$ is torsion-free and concentrated in degree $n$.  If a group $\gG$ acts freely and cocompactly on a duality space, then it is a duality group.

A finite simplicial complex $L$ has \emph{punctured (co)homology concentrated in degree $m$} if, for each closed simplex of $L$ (including the empty simplex), the reduced cohomology, $\ol{H}^*(L-\gs)$, is torsion-free and concentrated in degree $m$ (cf Davis--Meier \cite[\S6]{dm}).  We write \emph{$L$ is $PH^m$} as a shorthand for this condition.  For example, any triangulation of $S^m$ is $PH^m$.

The following proposition is proved in Davis et al \cite{ddjmo} (also see Davis--Meier \cite{dm,dm-erratum} and Davis et al \cite{ddjo2}).

\begin{proposition}\label{p:d}
\textup{(Davis et al \cite[Cor.~9.5]{ddjmo})}
Suppose $\cac$ is a building of type $(W,S)$.  Then its standard realization, $\cu(\cac,K)$, is an $n$-dimensional duality space if and only if the nerve, $L(W,S)$,  is $PH^{n-1}$.
\end{proposition}

\section{Constructing examples using covering spaces}\label{s:ex}
\subsection{The main construction}
Next we recall a construction of Davis \cite[\S 2]{d87}.  Suppose $L$ is a simplicial complex with vertex set $S$,  $(W',S')$ is a Coxeter system and  $f:S\to S'$ is a function defining a simplicial map (also denoted by $f$) from $L$ to $L(W',S')$.  Also suppose that the restriction of $f$ is each simplex is injective.  (Most of the time, $f:S\to S'$ will be a bijection and $W'$ will be a finite Coxeter group.)  Let $(m'(s',t'))$ be the Coxeter matrix for $(W',S')$.  Define an $(S\times S)$ Coxeter matrix $(m(s,t))$ by
\begin{equation}\label{e:coxm}
m(s,t):=
\begin{cases}
1,	&\text{if $s=t$,}\\
m'(f(s),f(t)),	&\text{if $\{s,t\}\in \edge(L)$,}\\
\infty,	&\text{otherwise}.
\end{cases}
\end{equation}
Let $(W,S)$ be the corresponding Coxeter system.  The map $f:S\to S'$ of generating sets extends to a homomorphism $\gf_f:W\to W'$. Let $\pi:=\Ker (\gf_f)$.  The space $K(L)$ has a mirror structure over $S$ and so, by  using $f$, we also get one over $S'$, defined by
\begin{equation*}
K(L)_{s'}:= \bigcup_{s\in f^\minus(s')} K(L)_s.
\end{equation*}
As in \S\ref{ss:gr}, we have the $W$-space, $\cu(W,K(L))$, and the $W'$-space, $\cu (W', K(L))$.  The space $\cu(W', K(L))$ is connected if and only if $f:S\to S'$ is surjective (cf Davis \cite[Remark 2.2]{d87}).  Henceforth, assume this.  Then $\gf_f:W\to W'$ is an epimorphism and it induces  $\gf_f$-equivariant map $q:\cu(W,K(L))\to \cu(W', K(L))$.  It is easy to see that $q$ is a covering projection and that $W$ is the group of all lifts of the $W'$-action.  By Davis \cite[Thm. 9.1.3, p. 166]{dbook}, $\cu(W,K(L))$ is simply connected.  Hence,
\[
\pi=\pi_1(\cu(W', K(L)).
\]
\begin{Remark}
$\cu(W, K(L))$ need not be contractible.  A necessary and sufficient condition for this to be true is that for every spherical subset  $T\subset S$, $T$ is the vertex set of a simplex in $L$, i.e, $L=L(W,S)$.
\end{Remark}

Next we want to carry out the same construction for buildings.  Suppose we are given the same data as above, as well as, a building $\cac'$ of type $(W',S')$.  Let $v_\emptyset$ be the central vertex of $K(L)$ corresponding to $\emptyset$.   For each $C'\in \cac'$, let $v_{C'}$ denote the corresponding central vertex of $\cu(\cac',K(L))$.  We can identify $\cac'$ with the set $\cent(\cu(\cac', K(L))):=\{v_{C'}\}_{C'\in \cac'}$ of central vertices.  By analogy with the case of Coxeter groups, consider the universal cover $p:\tcu\to\cu(\cac', K(L))$.  Define
\begin{equation}\label{e:cac}
\cac:=\cent(\tcu):=p^\minus (\cent(\cu(\cac', K(L))).
\end{equation}
Let $\tilde{\pi}$ denote the fundamental group of $\cu(\cac',K(L))$.  There is a free action of $\tilde\pi$ on $\cac$ such that the quotient set is identified with $\cac'$.  Since $K(L)$ is simply connected (it is a cone), $p^\minus (K(L))$ is isomorphic to $\tilde\pi \times K(L)$.  For each $C\in \cac$, let $\K$ denote the component of $p^\minus (K(L))$ which contains the vertex $C$.  We shall say that $\K$ is a \emph{chamber} of $\tcu$.  Note that $p\vert_{\K}$ maps $\K$ homeomorphically onto $K(L)$.

We want to show that $\cac$ has the structure of a building of type $(W,S)$ and that $\tcu\cong\cu(\cac, K(L))$.  To this end, we first show $\cac$ has the structure of a chamber system over $S$.   Chambers $C,D\in \cac$ are \emph{$s$-adjacent} if $\K\cap D_{\langle K\rangle}$ is nonempty and projects homeomorphically onto $K(L)_{f(s)}$.   A gallery in $\cac$ ($=\cent (\tcu)$) is \emph{reduced} if its type $\bs$ is a reduced expression for $w(\bs)$.

\begin{lemma}\label{l:gallery}
Let $\cac=\cent(\tcu)$ be defined by \eqref{e:cac}.  Suppose $(C_0, \dots, C_k)$ is a gallery in $\tcu$ of type $\bs$ and that $\bprime$ is another word in $S$ which can be obtained from $\bs$ by flip moves. Then there is a gallery in $\tcu$ of type $\bprime$ with the same endpoints.  
\end{lemma}

\begin{proof}
Suppose $T\in\cs(W,S)$ is a spherical subset.  Then the restriction of  $p$ to a $T$-residue in $\cac$ is a bijection to   (the set of central vertices in) a $T$-residue of $\cac'$.  Hence, every $T$-residue in $\cac$ is a spherical building.  It suffices to show that we can accomplish a single flip move on $\bs$ without changing the endpoints.  So, suppose $\bs=\Prod(s,t;m(s,t))$, where $m(s,t)<\infty$.  Then the entire gallery from lies in a single $\{s,t\}$-residue.  Since this residue is a building and $\Prod(s,t;m(s,t))$ is a reduced expression, there is a gallery of type $\Prod(t,s;m(s,t))$ with the same endpoints. In the general case we apply this step to the subgallery in which the flip is supposed to take place,  proving the lemma.
\end{proof}

Tits' solution to the Word Problem for Coxeter groups gives the following. 
\begin{corollary}\label{c:gallery}
Suppose $(C_0, \dots, C_k)$ is a gallery in $\tcu$ of type $\bs$ with $\bs$ a reduced expression for an element of $W$.  Let $\bprime$ be another reduced expression for $w(\bs)$.  Then there is a gallery of type $\bprime$ from $C_0$ to $C_k$.
\end{corollary}

\begin{lemma}\label{l:mingallery}
Suppose $\bC=(C_0,\dots, C_k)$ is a gallery of type $\bs$ from $C_0$ to $C_k$.  If $\bC$ has minimum length, then $\bs$ is a reduced expression.
\end{lemma}

\begin{proof}
If $\bs$ is not reduced, then, by Tits' solution to the Word Problem, $\bs$ can be changed by a sequence of flips to $\bprime=(s'_1,\dots,s'_k)$ so that $\bprime$ contains a subword of length two of the form $(s'_j,s'_{j+1})=(s,s)$.  By Lemma~\ref{l:gallery}, there is a gallery $\bC'=(C'_0,\dots, C'_k)$ with the same endpoints as $\bC$.  Let $(C'_{j-1}, C'_j, C'_{j+1})$ be the subgallery corresponding to the subword $(s,s)$.  Then $C'_{j-1}, C'_j$ and $C'_{j+1}$ are all $s$-equivalent.  If $C'_{j+1}=C'_{j-1}$, then we can shorten the gallery by deleting $C'_{j-1}$ and $C'_j$.  If not, we can shorten it by deleting $C'_j$.  In either case we contradict the assumption that $\bC$ has minimum length.
\end{proof}

Next, we want to define the Weyl distance $\gd:\cac\times \cac \to W$.  Fix a chamber $C'\in \cac'$ and let $\gr_{C'}:\cac'\to W'$ be the function $D'\mapsto \gd'(C',D')$.  Since $\gr_{C'}$ maps $s'$-equivalent chambers to $s'$-equivalent elements of $W'$, it induces a map of spaces $r_{C'}:\cu(\cac', K(L))\to \cu(W',K(L))$.  Choose a chamber $C\in \cac$ lying over $C'$.  We want to construct  a map $r_C:\tcu \to \cu(W,K(L))$ making the following diagram commute:
\[
\begin{CD}
\tcu@>r_{C}>> \cu(W,K(L))\\
@Vp VV	@VVqV\\
\cu(\cac',K(L)))@>r_{C'} >>\cu(W',K(L))
\end{CD}
\]
Since $q$ is a covering projection and $\cu(W,K(L))$ is simply connected, the  standard lifting theorem in covering space theory implies that  $r_{C'}\circ p: \tcu \to \cu(W',K(L))$ has a unique lift $r_C$, as in the diagram, and  taking the central vertex $C$ to $1$.  By restricting $r_C$ to the subset  $\cac\subset\tcu$ we get a function $\gr_C:\cac\to W$.  ($\cac$ is the set of central vertices in $\tcu$ and $W$ is the set of central vertices in $\cu(W,K(L))$.)  Define $\gd:\cac\to W$ by
\begin{equation}\label{e:Wdistance}
\gd(C,D):=\gr_C(D).
\end{equation}

\begin{Remark}
The reason this procedure works comes down to a fact, which was used in the proof of Lemma~\ref{l:gallery}: every spherical residue in 
the chamber system 
$\cac$ is a spherical building.  This together with simple connectivity of the standard realization implies that $\cac$ is a building (cf Tits \cite{tits86}). 
\end{Remark}

\begin{lemma}\label{l:Wdistance}
Suppose $\bC=(C_0,\dots, C_k)$ is a reduced gallery in $\cac$ of type $\bs$.  Then $\gd(C_0,C_k)=w(\bs)$.
\end{lemma}

\begin{proof}
Let $\gamma:[0,k]\to \tcu$ be the edge path in the dual $1$-skeleton of $\tcu$ corresponding to $\bC$.  Thus, $\gamma$ crosses the central vertices $C_0, C_1\dots, C_k$ in succession.  Then $r_{C_0}\circ \gamma$ is an edge path in the dual $1$-skeleton of $\cu(W,K(L))$.  Let $(w_0, w_1, \dots, w_k)$ be the corresponding gallery in $W$.  By construction the type of this gallery is $\bs$ (the type of $\bC$).  Since $r_{c_0}(C_0)=1$, the initial vertex $w_0$ is $1$.  Thus, $\gd(C_0, C_k)=r_{C_0}(C_k)=w_k=w(\bs)$.
\end{proof}

\begin{theorem}\label{t:cac}
$\gd:\cac\times\cac\to W$ is a Weyl distance function; hence, $\cac$ is a building of type $(W,S)$.  Moreover, $\tcu$ can be identified with $\cu(\cac,K(L))$.
\end{theorem}

\begin{proof}
We must check conditions \textbf{(WD1)}, \textbf{(WD2)}, \textbf{(WD3)} of Definition~\ref{d:bldg}.  The first condition is immediate.  To check \textbf{(WD2)} suppose $\gd(C,D)=w$ and $C'\in \cac$ satisfies $\gd(C'C)=s$.   By definition of $\gd$ there is a gallery $(C_0, \dots, C_k)$ of type $\bs$ from $C$ to $D$, where $\bs=(s_1,\dots, s_k)$ is a reduced expression for $w$. There are two cases to consider.

\textbf{Case 1:} $l(sw)>l(w)$.

Then $\bprime=(s,s_1,\dots, s_k)$ is a reduced expression for $sw$ and $(C',C_0,\dots, C_k)$ is a reduced gallery of type $\bprime$.  So, by Lemma~\ref{l:Wdistance}, $\gd(C',C_k)=w(\bprime)=sw$.
\textbf{Case 2:} $l(sw)<l(w)$.

Then there is a reduced expression $\bprime$ for $w$ which begins with $s$.  By Corollary~\ref{c:gallery}, there is a gallery of type $\bprime$ from $C$ to $D$.  The resulting  gallery $(C', C, C_1,\dots, D)$ has type beginning $(s,s)$.  Hence, the chambers $C'$, $C$ and $C_1$ are $s$-equivalent.  If $C_1=C'$, then $(C_1, C_2,\dots)$ is a reduced gallery from $C'$ to $D$; hence, $\gd(C',D)=sw$.  If $C_1\neq C$, then $(C', C_1, C_2,\dots)$ has type $\bprime$ and is a gallery from $C'$ to $D$.  Hence, $\gd((C',D)=w(\bprime)=w$.

To check \textbf{(WD3)}, let $s$ be an arbitrary element of $S$.  Suppose we have a reduced gallery  $(C_0, \dots, C_k)$ from $C$ to $D$ of type $(s_1, \dots, s_k)$.  The argument divides into the same two cases as above.  In Case 1, $(C', C_0, C_1,\dots, C_k)$ is a reduced gallery; hence, $\gd(C',C_k)=sw$.  In Case 2, by using Corollary~\ref{c:gallery} and Lemma~\ref{l:Wdistance},  we can assume that $s_1=s$ (ie, the type of the gallery begins with $s$).  Then choosing $C'=C_1$ we get a reduced gallery $(C_1, \dots, C_k)$ from $C'$ to $D$; hence, $\gd(C',D)=sw$.
\end{proof}

\begin{examples}\label{ex:fold}
Suppose $W'$ is a finite Coxeter group of rank $n+1$ with generators indexed by $S'=\{0,1,\dots,n\}$.  Then  $L'$($=L(W',S')$) is the simplex on $\{0,1,\dots,n\}$.  Suppose $L$ is an $n$-dimensional simplicial complex which admits a ``folding map'' $f:L\to L'$ (ie, the restriction of $f$ to each simplex is injective).  Then we can proceed as above.  In particular:

\begin{enumerate1}
\item
Suppose $L$ is the barycentric subdivision of an $n$-dimensional cell complex.  Then there is a folding map $f:S\to S'$ which sends the barycenter of a cell $F$ to $\dim F$.  If $L$ is the barycentric subdivision of the boundary complex of an $(n+1)$-dimensional polytope, then $K(L)$ is a simple polytope.  For example, if $L$ is the boundary of a $(n+1)$-simplex, then $K(L)$ is an $(n+1)$-dimensional permutohedron (cf Davis \cite[Remark 6.6]{d87}).  In the building $\cu(\cac,K(L))$ the link of each vertex of $K(L)$ is the spherical building $\cu(\cac', \gD^n)$, where $\gD^n$ is the dual simplex to $L'$.
\item
Suppose $L$ is the boundary complex of an $(n+1)$-dimensional octahedron (see \S\ref{ss:square} below).  Its vertex set, $S$, can be identified with $\{0,\dots, n\}\times \{\pm 1\}$ and the folding map is induced by projection on the first factor, $\{0,\dots, n\}\times \{\pm 1\}\to \{0,\dots, n\}$.  In this case, $K(L)$ is an $(n+1)$-cube.  As in 1), the link of each vertex of $K(L)$ in $\cu(\cac,K(L))$ is the spherical building $\cu(\cac', \gD^n)$.
\end{enumerate1}
\end{examples}

\begin{definition}\label{d:flag}
A simplicial complex $J$ is a \emph{flag complex} if  any nonempty, finite set of vertices, which are pairwise connected by edges, spans a simplex of $J$.
\end{definition}

A flag complex $J$ is ``determined by its $1$-skeleton'' in that it is the smallest full subcomplex of the full simplex on $\vertex (J)$ which contains the $1$-skeleton $J^1$.

\begin{theorem}\label{t:flag}
If in the above construction $L$ is a flag complex, then $L=L(W,S)$ (and hence, $\cu(\cac,K(L))$ is the standard realization of $\cac$).
\end{theorem}

\begin{proof}
Suppose $T$ is a spherical subset of $S$.  Since this implies that $m(s,t)<\infty$ for all $s,t\in T$, we see that the elements of $T$ are pairwise connected by edges in $L$.  Since $L$ is a flag complex, this means that $T$ is the vertx set of a simplex in $L$.  Hence, $L(W,S)=L$.
\end{proof}

\begin{Remark}
When $W'$ is finite the condition that $L$ be a flag complex is also necessary for $L$ to equal $L(W,S)$.
\end{Remark}

\begin{example}\label{ex:random}
Given $(W',S')$, choose a set of edges $E\subset \edge(L(W',S'))$.  Let $f:S\to S'$ be any bijection.  As in \eqref{e:coxm}, define a Coxeter matrix $(m(s,t))$ by
\[
m(s,t)=
\begin{cases}
m'(f(s),f(t)),	&\text{if $\{s,t\}\notin E$,}\\
\infty,	&\text{if $\{s,t\}\in E$.}
\end{cases}
\]
Let $(W,S)$ be the resulting Coxeter system and $L=L(W,S)$. So, if $\cac'$ is a building of type $(W',S')$, we get a building $\cac$ of type $(W,S)$.

Here is slightly different spin on the same example.  Suppose $L$ is an arbitrary finite flag complex with vertex set $S$.  Then $L$ is a subcomplex of the full simplex $\gD$ on $S$.  Let $(W',S)$ be an arbitrary spherical Coxeter system.  Then $L(W',S)=\gD$.  If $E$ is the set of edges of $\gD$ which are not in $L$, then the construction in the previous paragraph yields $(W,S)$ with $L(W,S)=L$.  For $L\neq \gD$, we could have reached the same conclusion by choosing $(W',S)$ to be any Coxeter system with nerve $\partial \gD$.  For example, $(W',S')$ could be an irreducible affine Coxeter system.

Suppose the flag complex $L$ is $PH^{n-1}$ (defined in \S \ref{ss:duality}).  Then, by Proposition~\ref{p:d}, $\cu(\cac,K(L))$ is an $n$-dimensional duality space and, by Corollary~\ref{c:lattice} below, $\Aut (\cac)$ is an $n$-dimensional duality group.
\end{example}

\subsection{Automorphism groups}
\begin{proposition}\label{p:transitive}
Suppose $G'$ is a group of automorphisms of $\cac'$.  Then $\cac$  has a group of automorphisms $G$ with $\cac/G\cong \cac'/G'$ and $\cu(\cac,K(L))/G= \cu(\cac', K(L))/G'$.  
\end{proposition}

\begin{proof}
By Theorem~\ref{t:cac}, 
$\cu(\cac,K(L))\to \cu(\cac', K(L))$ is the universal covering space. Let $G$ be the group of all lifts of the $G'$-action.  The asserted isomorphisms are easily established .
\end{proof}

\begin{corollary}\label{c:transitive}
If $\cac'$ admits a chamber-transitive group of automorphisms,  then so does $\cac$.     More generally, if the action of a group $G'$ of automorphisms of $\cac'$ has only finitely many orbits, then $G$ acts on $\cac$ with the same number of orbits. 
\end{corollary}
 
\begin{corollary}\label{c:lattice}
If $\cac'$ is a spherical building, then $\Aut(\cac)$ admits a torsion-free uniform lattice. 
\end{corollary}

\begin{proof}
We claim that the discrete subgroup $\tilde{\pi}\subset \Aut(\cac)$ is such a lattice (where $\tilde{\pi}:=\pi_1(\cu(\cac',K)$).  It is uniform, since $\cac/\tilde{\pi}=\cac'$, which is a finite set.  It is torsion-free, since it acts freely on the finite dimensional, contractible space $\cu(\cac,K)$.
\end{proof}

\begin{Question}
When is $\Aut(\cac)$ strongly transitive?
\end{Question}

\section{Partial products of Coxeter groups and buildings}\label{s:ppc}
\subsection{General description}
We are given as data:
\begin{itemize}
\item
Coxeter systems $(W_i, S_i)$, for $i=0,1,\dots, p$,
\item
for each pair $(i,j)$ with $0\le i,j\le p$ and $i\neq j$, a subset $\car_{ij}\subset S_i\times S_j$ such that $\car_{ji}$ is obtained from $\car_{ij}$ simply by interchanging the factors (ie, $(s_j,s_i)\in \car_{ji}$ $\iff$ $(s_i,s_j)\in \car_{ij}$).
\end{itemize}

From this data we  define a new Coxeter system $(W,S)$.  First, $S$ is the disjoint union $S:=S_0\coprod \dots \coprod S_p$.  For $s_i,t_i\in S_i$, let $m_i(s_i,t_i)$ denote the order of $s_it_i$.   Define an $(S\times S)$ symmetric matrix $(m(s,t))$ as follows:  given $s\in S_i$, $t\in S_j$,
\[
m(s,t)=
\begin{cases}
m_i(s,t),	&\text{if $i=j$,}\\
\infty,	&\text{if $i\neq j$ and $(s,t)\in \car_{ij}$,}\\
2,	&\text{if $i\neq j$ and $(s,t) \notin \car_{ij}$.}
\end{cases}
\]
$W$ is the Coxeter group associated to the Coxeter matrix $(m(s,t))$.
In other words, the $\car_{ij}$ specify the edges from a vertex in $S_i$ to one in $S_j$ which are to be deleted from the nerve of $(W_0\times\cdots\times W_p, S)$.

\begin{examples}\label{ex:first}
(1) Suppose $\car_{ij}=\emptyset$ for all $i\neq j$.  Then $W$ is the direct product $W_0\times W_1\times\cdots \times W_p$.

(2) Suppose $\car_{ij}=S_i\times S_j$ for all $i\neq j$.  Then $W$ is the free product $W_0*W_1*\cdots *W_p$.

(3) Suppose $\gO$ is a simplicial graph with vertex set $\{0,1,\dots, p\}$ and 
\[
\car_{ij}=
\begin{cases}
\emptyset,	&\text{if $\{i,j\} \in \edge(\gO)$,}\\
S_i\times S_j,			&\text{if $\{i,j\}\notin \edge(\gO)$.}
\end{cases}
\]
Then $W$ is the ``graph product'' of the $W_i$, denoted by $\prod_\gO W_i$.  (See Januszkiewicz-- \'Swi\c atkowski \cite{js} for discussion of graph products.)  As a further specialization, if each $W_i$ is cyclic of order two, then $W=\prod_\gO \zz/2$ is the right-angled Coxeter group associated to $\gO$.  If each $W_i$ is the infinite dihedral group, $D_\infty$, then $W=\prod_\gO D_\infty$ contains the right-angled Artin group associated to $\gO$ as a subgroup of index $2^{p+1}$ (cf Davis--Januszkiewicz \cite{dj00}).
\end{examples}

\begin{example}\label{ex:square}
(\emph{Square products})\ \   In this example suppose the sets $S_i$ are mutually bijective and that we are given 
a compatible family of bijections, $\{\theta_{ij}:S_i\to S_j\}_{0\le i,j\le p}$.
(By a ``compatible family'' we mean that $\theta_{ii}=1$ and $\theta_{ jk}\circ \theta_{ij}=\theta_{ik}$, for all $i$, $j$, $k$.)  The $\car_{ij}$ are required to satisfy the following:  the bijection $\theta_{i0}\times \theta_{j0}:S_i\times S_j\to S_0\times S_0$ takes $\car_{ij}$ onto the diagonal of $S_0\times S_0$ (or equivalently, $(\theta_{ij}\times id)(\car_{ij})$ is the diagonal of $S_j\times S_j$). Then for $s\in S_i$, $t\in S_j$,
\[
m(s,t)=
\begin{cases}
m_i(s,t),	&\text{if $i=j$,}\\
\infty,		&\text{if $i\neq j$ and $\theta_{ij}(s)=t$,}\\
2,	&\text{otherwise.}
\end{cases}
\]
In this case call $W$ the \emph{square product} of the $W_i$ and write $W=W_0\,\square\cdots \square W_p$. (In the case where there are only two factors, $W_0$, $W_1$, and both are finite, these examples were discussed in Davis \cite[Example 4.3]{d87}.)
\end{example}

\textbf{Spherical subsets of $\boldsymbol{S}$}\  \  Let $\cs_i:=\cs(W_i,S_i)$ be the poset of spherical subsets of $S_i$ and let $\cs:=\cs(W,S)$ be the poset of spherical subsets of $S$.  Any $T\subset S$ can be decomposed as $T=T_0\coprod \cdots \coprod T_p$, where $T_i\subset S_i$.  Clearly, $T$ is spherical if and only if the following two conditions hold:
\begin{itemize}
\item
$T_i\in \cs_i$, for $i=0,\dots,p$, and
\item
$(\theta_{i0}(s_i),\theta_{j0}(s_j))\notin \car_{ij}$, for all $s_i\in T_i$, $s_j\in T_j$, $i\neq j$.
\end{itemize}
Put 
\begin{equation}\label{e:W'}
W':=W_0\times \cdots \times W_p
\end{equation}
and note that 
\(
\cs(W',S)=\cs_0\times \cdots\times \cs_p
\). As in \S 2, the tautological bijection $f:S\to S_0\coprod\cdots\coprod S_p$ extends to an epimorpism $\gf_f:W\to W'$ and induces an embedding $\cs\to \cs(W',S)$ defined by $T\mapsto (T_0,\dots, T_p)$ (where $T=T_0\coprod\cdots \coprod T_p$).
Similarly, if  $L:=L(W,S)$ and $L_i:=L(W_i,S_i)$ denote the nerves of their respective Coxeter systems, we get a simplicial embedding $f:L\to L(W',S)$.  (Note that $L(W',S)$ is the join, $L_0*\cdots *L_p$.)

\textbf{Buildings}\ \   For $0\le i\le p$, suppose $\cac_i$ is a building of type $(W_i,S_i)$.  Then $\cac_0\times\cdots\times  \cac_p$ is a building of type $(W', S)$.  As in \S 2 this lead to a new building $\cac$ and a covering space, $\cu(\cac,K)\to \cu(\cac',K)$ where $K:=K(L)$ is the fundamental chamber for $(W,S)$.

If $W$ is as in Examples~\ref{ex:first} part 1), 2) or 3), then  $\cac$ is, respectively,  the \emph{product}, \emph{free product} or \emph{graph product} of the $\cac_i$.  In the case of the graph product, write $\cac=\prod_\gO \cac_i$. 
If $W$ is as in Example~\ref{ex:square}, then $\cac=\cac_0\,\square\cdots\square\, \cac_p$ is the \emph{square product} of the $\cac_i$.  

We have the following corollary to Corollary~\ref{c:lattice}.

\begin{corollary}\label{c:lattice2}
Suppose each $\cac_i$ is a spherical building (of finite thickness) and $\gO$ is a simplicial graph.  If  $\gO$ is not the $1$-skeleton of a simplex, then the graph product $\prod_{\gO} \cac_i$ is an infinite building whose automorphism group admits a uniform lattice.  Similarly, if $p\ge 1$, then $\cac_0\,\square\cdots \square\, \cac_p$ is an infinite building whose automorphism group admits a uniform lattice.
\end{corollary}

\subsection{Description of $L$ for $\square$-products}\label{ss:square}
Suppose $V$ is a set with $n$ elements and $\gD(V)$ is the simplex on $V$.  Define a mirror structure over $V$ on $\gD(V)$ by letting $\gD(V)_v$ be the codimension one face opposite to the vertex $v$.  Let $(\zz/2)^V$ be the direct product of $V$ cyclic groups of order two.  Identify $V$ with the standard basis of $(\zz/2)^V$ (so that $((\zz/2))^V,V)$ is a Coxeter system). Then $\cu((\zz/2)^V,\gD(V))$ is a triangulation of  $S^{n-1}$, specifically, the boundary complex of an $n$-dimensional octahedron.

For any integer $p\ge 1$, the set $\{0,\dots, p\}^V$ is a building of type $((\zz/2)^V,V)$ and thickness $p$.  Its realiztion $\cu(\{0,\dots,p\}^V,\gD(V))$ is an $n$-dimensional spherical building.  The space $\cu(\{0,\dots,p\}^V,\gD(V))$ is a simplicial complex which we will denote $O(p,V)$.  The simplices of $O(p,V)$ can be descibed as follows. 
Consider pairs of the form $g\times \gD(U)$ where $U\subset V$ and where $g=(g_v)_{v\in V}$ is a $V$-tuple in $\{0,\dots, p\}$.  Two such pairs $g\times \gD(U)$ and $g'\times \gD(U')$ are equivalent if and only $U=U'$ and $g_v=g'_v$ for all $v\in U$.  A simplex in $O(p,V)$ is an equivalence class.

The complex $O(p,V)$ can also be described as the $n$-fold join of the set $\{0,\dots,p\}$ with itself. 

\begin{lemma}\label{l:ph}
$O(p,V)$ is $PH^{n-1}$.
\end{lemma}

\begin{proof}
The vertex set of an $(m-1)$-simplex $\gs$ in $O(p,V)$ is a $U$-tuple $g\in \{0,\dots,p\}^U$, where $m=\Card(U)$.  The complement of $\gs$ in $O(p,V)$ deformation retracts onto the complement of the open star of $\gs$.  The complement of this open star can be identified with the join of $m$ sets of $p$ elements and $n-m$ sets with $p+1$
It is well-known that the join of $n$ finite sets, each of which has at least $2$ elements, is homotopy equivalent to a wedge of $(n-1)$-spheres.  If at least one of finite, nonempty sets has exactly one point, then the join is contractible.  The lemma follows. 
\end{proof}

As before, $(W_i,S_i)$ are Coxeter systems, $\theta_{ij}:S_i\to S_j$ are compatible bijections, $L_i=L(W_i,S_i)$ and $\cs_i=\cs(W_i,S_i)$.  Also, $S:=S_0\coprod\cdots \coprod S_p$ and $W:=W_0\,\square\cdots\square\,W_p$
Put $V:=S_0$.  Using the bijections $\theta_{i0}:S_i\to S_0$  identify each $S_i$ with $V$.
We are going to define an embedding of posets $\cs\hookrightarrow \cs(O(p,V))$.  The image will be denoted $\cs_0\,\square\cdots \square\,\cs_p$.  An element of $\cs$ has the form $(T_0, \dots,T_p)$, where $T_i\subset S_i$.  Use the $\theta_{i0}$ to transport the $T_i$ to subsets $U_i$ of $V$, ie, put $U_i:=\theta_{i0}(T_i)$.  The fact that each $\car_{ij}$ is the diagonal of $V\times V$ means that the $U_i$ are pairwise disjoint subsets of $V$.  Put $U:=U_0\cup \cdots \cup U_p$.  Since the $U_i$ are disjoint, any $u\in U$ belongs to exactly one of them, call it $U_{g_u}$.  This defines a $U$-tuple $g\in \{0,\dots, p\}^U$.  Extend this arbitrarily to a $V$-tuple in $\{0,\dots, p\}$, also denoted by $g$. (In other words, $g_v$ is completely arbitrary for $v\in V-U$.)  The embedding  $\cs\hookrightarrow \cs(O(p,V))$ is then defined by $(T_0,\dots, T_p)\mapsto (U,g)$.  The corresponding subcomplex of $O(p,V)$ is denoted $L_0\,\square\cdots\square\,L_p$.  The above discussion proves the following.

\begin{proposition}\label{p:square}
With notation as above, 
\[
L(W,S)\cong L_0\,\square\cdots\square\,L_p\subset O(p,V).
\]
\end{proposition}

\begin{Remark}
The above description of $L=L_0\,\square\cdots\square\,L_p$ gives the following estimate for its dimension:
\[
\max_i\{\dim L_i\}\le\dim L\le n-1.
\]
\end{Remark}

\begin{corollary}\label{c:square}
Suppose each $W_i$ is a finite Coxeter group of rank $n$ and $V$ is a set with $n$ elements.  Then $L(W,S)\cong O(p,V)$. 
\end{corollary}

Combining this with Lemma~\ref{l:ph} and Proposition~\ref{p:d} we get the following.

\begin{corollary}
Suppose each $W_i$ is a finite Coxeter group of rank $n$ so that each $\cac_i$ is a spherical building. Then the standard realization of $\cac_0\,\square\cdots\square\,\cac_p$ is an $n$-dimensional duality space and 
$\Aut(\cac_0\,\square\cdots\square\,\cac_p)$ is a virtual $n$-dimensional duality group.
\end{corollary}

\begin{example}\label{ex:cube}
When $p=1$, $O(1,V)$ is the boundary complex of an $n$-dimensional octahedron.  So, in this case, $K$ ($=K(L)$)  is combinatorially isomorphic to an $n$-cube.  Suppose $W_0$, $W_1$ are finite Coxeter groups of rank $n$.  Then $K$ is the $n$-cube and in the standard realization, $\cu(W_0\,\square\, W_1, K)$, the isotropy group at one vertex of $K$ is $W_0$ while at the opposite vertex it is $W_1$, cf Davis \cite[Example~4.3]{d87}.  It follows from this that $W_0K$ is combinatorially isomorphic to a Coxeter cell of type $(W_0,S_0)$.  For example, if $W_0$ is the symmetric group of rank $n+1$, then $W_0K$ is a permutohedron, cf \cite[\S6]{d87}.  Next suppose that, for $i=0,1$, $\cac_i$ is a spherical building of type $(W_i,S_i)$.  Then in the standard realization, $\cu(\cac_0\, \square\, \cac_1)$, the link at one vertex is the spherical realization $\cu(\cac_0, \gD)$ while at the other it is $\cu(\cac_1, \gD)$.  At the other vertices, the links are joins of the form $\cu(\cac'_0,\gD_0)*\cu(\cac'_1,\gD_1)$, where for $i=0,1$, $\cac_i$ is a building of type $(W_{T_i}, T_i)$ for $T_i\subset S_i$ and $\dim \gD_i=\Card (T_i)-1$.  If we choose $\cac_0$ to be the thin building $W_0$, then $W_0\, \square\, \cac_1$ can be regarded as the building with fundamental chamber $W_0K$ obtained, as in Example~\ref{ex:fold}, by using the natural folding map from the Coxeter complex, $\cu(W_0, \gD)$, to $\gD$.
\end{example}

\section{Questions related to Kac-Moody groups}
Suppose $(W,S)$ is such that $m(s,t)\in \{2,3,4,6,\infty\}$ for all $s\neq t$ and that $q$ is a prime power.  As in Tits \cite{tits87} one can choose an $(S\times S)$ ``Cartan matrix'' compatible with $(W,S)$ and from this construct a Kac-Moody Lie algebra and a Kac-Moody group $G$ over the field $\ff_q$ of order $q$.  Associated to $G$ there is a Tits system $(G,B,N, T)$ and a building $\cac_G$ with set of chambers $G/B$.  Call $\cac_G$ a \emph{Kac-Moody building}.

Now suppose as in  \S \ref{s:ex} or \S\ref{s:ppc} that $\cac'$ is a spherical building of type $(W',S')$ associated to an algebraic group $G'$ over a finite field $\ff_q$.  Then $(W,S)$ satisfies the crystallographic condition (that $m(s,t)\in \{2,3,4,6,\infty\}$).  Let $\cac$ be the building constructed in \S \ref{s:ex}.  The question arises:  is $\cac$ isomorphic to a Kac-Moody building $\cac_G$.  Undoubtedly, it often is.  For example, if $(W,S)$ is is a right-angled and $\cac$ is of thickness $q$ (so that $G'$ is the $S'$-fold product, $\prod PGL(2,\ff_q)$), then  $\cac\cong \cac_G$ (since any two regular right-angled buildings of type $(W,S)$ and thickness $q$ are isomorphic)  However, even when $\cac\cong \cac_G$, the group $\wt{G}\subset \Aut(\cac)$, obtained by lifting the $G'$-action to $\cac$, will not be closely related to the Kac-Moody group $G$.  One reason is that if $(W,S)$ is irreducible, infinite and nonaffine, then, by a result of Caprace-Remy \cite{cr}, $G$ is a simple group.  On the other hand, there is an epimorphism $\wt{G}\to G'$.  Hence, $G$ is contained in $\Ker (\wt{G}\to G')$. 

\bibliographystyle{gtart.bst}

\end{document}